\documentclass[10pt,aps,prb, preprint,reqno]{amsart}
\usepackage{amsmath}
\usepackage{amssymb}
\usepackage{hyperref}

\newtheorem{theorem}{Theorem}[section]
\newtheorem{remark}{Remark}[section]
\newtheorem{lemma}[theorem]{Lemma}
\newtheorem{proposition}[theorem]{Proposition}
\newtheorem{corollary}[theorem]{Corollary}

\begin{document}
\title[Regularity Criterion of 3D MHD system]{Regularity criteria of the three-dimensional MHD system involving one velocity and one vorticity component}
\author{Kazuo Yamazaki}  
\date{}
\maketitle

\begin{abstract}
We obtain a regularity criteria of the solution to the three-dimensional magnetohydrodynamics system to remain smooth for all time involving only one velocity and one vorticity component. Moreover, the norm in space and time with which we impose our criteria for the vorticity component is at the scaling invariant level. The proof requires a new decomposition of the four non-linear terms making use of a new identity due to the divergence-free conditions of the velocity and the magnetic vector fields.  

\vspace{5mm}

\end{abstract}
\footnote{2000MSC : 35B65, 35Q35, 35Q86}
\footnote{Address: Department of Mathematics, Washington State University, Pullman, WA 99164-3113; USA; phone: 1-509-335-9812; fax: (509) 335-1188; e-mail: kyamazaki@math.wsu.edu.	
	}

\section{Introduction}

We study the following magnetohydrodynamics (MHD) system: 

\begin{subequations}
\begin{align}
&\frac{du}{dt} + (u\cdot\nabla) u + \nabla \pi = \nu \Delta u + (b\cdot\nabla) b,\\
&\frac{db}{dt} + (u\cdot\nabla) b = \eta \Delta b + (b\cdot\nabla)u,\\
&\nabla \cdot u = \nabla\cdot b = 0, \hspace{3mm} (u,b)(x,0) = (u_{0}, b_{0})(x),
\end{align}
\end{subequations}
where $u: \mathbb{R}^{3} \times \mathbb{R}^{+} \mapsto \mathbb{R}^{3}, b: \mathbb{R}^{3} \times \mathbb{R}^{+} \mapsto \mathbb{R}^{3}, \pi: \mathbb{R}^{3} \times \mathbb{R}^{+} \mapsto \mathbb{R}$ represent the velocity, magnetic and pressure fields respectively. We have denoted by the parameters $\nu, \eta > 0$ the kinematic viscosity and magnetic diffusivity which is a reciprocal of the magnetic Reynolds number respectively; for simplicity we assume them to be one for the rest of the manuscript. The important study of investigating the motion of electrically conducting fluids can be traced back to the pioneering work in [1, 8]. Ever since then, the MHD system has found much applications in astrophysics, geophysics and plasma physics. We also remark that the system (1a)-(1c) at $b \equiv 0$ reduces to the Navier-Stokes equations (NSE).  

In order to be able to precisely describe previous work on the MHD system and the NSE, let us write components of $u, b$ by $u= (u_{1}, u_{2}, u_{3}), b =(b_{1}, b_{2}, b_{3})$ and the vorticity and current density respectively as follows: 

\begin{equation}
\omega = (\omega_{1}, \omega_{2}, \omega_{3}) \triangleq \nabla \times u, \hspace{5mm} j = (j_{1}, j_{2}, j_{3}) \triangleq \nabla \times b. 
\end{equation}
For brevity let us denote $\int_{\mathbb{R}^{3}} f dx$ by $\int f$, $\frac{d}{dt}$ by $\partial_{t}$ and  $\frac{d}{dx_{i}}$ by $\partial_{i}, i = 1, 2, 3$. With that, we also denote by $\nabla_{h} = (\partial_{1}, \partial_{2}, 0)$ and $\Delta_{h} = \sum_{k=1}^{2} \partial_{k}^{2}$. Finally, let us write $A \lesssim_{a,b} B$ when there exists a constant $c \geq 0$ of significant dependence only on $a, b$ such that $A \leq c B$, similarly $A \approx_{a,b} B$ in case $A = cB$. 

Due to the work in [19, 23], we know the global existence of a weak solution and local existence of the unique strong solution which in particular satisfies the following energy inequality: 
		
\begin{equation}
\left(\lVert u\rVert_{L^{2}}^{2} + \lVert b\rVert_{L^{2}}^{2}\right)(t) + \int_{0}^{t} \lVert \nabla u\rVert_{L^{2}}^{2} + \lVert \nabla b\rVert_{L^{2}}^{2} d\tau \leq \lVert u_{0}\rVert_{L^{2}}^{2} + \lVert b_{0}\rVert_{L^{2}}^{2}.
\end{equation}
However, the smoothness of such a weak solution to the MHD system or the NSE remains unknown and one of the best clues toward the resolution of the Navier-Stokes problem is the following result by the authors in [11, 24] and others  which showed that if a weak solution $u$ to the NSE on $[0,T]$ satisfies 
	
\begin{equation}	
\int_{0}^{T} \lVert u \rVert_{L^{p}}^{r}d\tau  < \infty, \hspace{3mm} \frac{3}{p} + \frac{2}{r} \leq 1, \hspace{3mm} p \in [3, \infty],
\end{equation}
then $u$ is smooth. We refer to [3] for other important result in this direction of research and also [13, 32] for extension to the MHD system. Let us also recall the following  important result from [4]: if $(u,b)$ is the unique local strong solution to (1a)-(1c) in $[0, T)$ and  

\begin{equation}
\int_{0}^{T} \lVert \omega\rVert_{L^{\infty}} + \lVert j\rVert_{L^{\infty}} d\tau < \infty, 
\end{equation}
then in fact, it remains a strong solution in $[0, T+\epsilon)$ for some $\epsilon > 0$. The case $j \equiv 0$ is the result from [2] on the NSE; we also note that the condition on $j$ may be eliminated even for the MHD system (e.g. [10]).   

We now focus on some component reduction results of such conditions that are of most relevance to our result  (cf. [5, 18, 22, 33] for the NSE, [6, 26, 27, 28, 29, 31] for the MHD system). Firstly, the authors in [17] (and also [5]) showed that for the NSE, upon the $\lVert \nabla_{h} u\rVert_{L^{2}}$-estimate, $u_{3}$ may be separated from the non-linear term as follows: 

\begin{equation}
\int (u\cdot\nabla) u \cdot \Delta_{h} u \lesssim \int \lvert u_{3} \rvert \lvert \nabla u\rvert \lvert \nabla\nabla_{h} u\rvert.
\end{equation}
With such a decomposition, the authors in [17] obtained the following component reduction result in comparison to (4) for the NSE:  

\begin{equation}
\int_{0}^{T} \lVert u_{3} \rVert_{L^{p}}^{r} d\tau < \infty, \hspace{3mm} \frac{3}{p} + \frac{2}{r} \leq \frac{5}{8}, \hspace{3mm} r \in [\frac{54}{23}, \frac{18}{5}]. 
\end{equation}
We refer to [14, 28] for extension to the MHD system that involves at least $u_{3}, b_{1}, b_{2}$. Moreover, importantly for our discussion, we remark that if $(u,b)(x,t)$ solves the MHD system, then so does 
$(u_{\lambda}, b_{\lambda})(x,t) \triangleq \lambda(u,b)(\lambda x, \lambda^{2} t)$ and such a solution is scaling-invariant under the norm $\int_{0}^{T} \lVert \cdot \rVert_{L^{p}}^{r} d\tau$ precisely when $\frac{3}{p} + \frac{2}{r} = 1$ and the upper bound of 1 in (4) has been compromised to $\frac{5}{8}$ in (7). To the best of the author's knowledge, a regularity criterion in terms of one velocity component for the MHD system even in non-scaling invariant level remains open, although e.g. in [26, 27, 31] a condition in terms of two components have been obtained. 

Another interesting result is a component reduction result from (5) in the case of the NSE. The authors in [7] showed in particular that a weak solution to the NSE is smooth if 

\begin{equation}
\int_{0}^{T} \lVert \omega_{2} \rVert_{L^{p}}^{r} + \lVert \omega_{3} \rVert_{L^{p}}^{r} d\tau < \infty, \hspace{3mm} \frac{3}{p} + \frac{2}{r} \leq 2, 1 < r < \infty  
\end{equation}
(cf. [16]). To the best of the author's knowledge, it is not known if a condition for the MHD system in terms of $\omega$ e.g. in [10] may be reduced to just $\omega_{2}, \omega_{3}$. The most similar result in this direction of research is the criterion of 

\begin{equation}
\int_{0}^{T} \lVert \partial_{3} u_{3} \rVert_{L^{p}}^{r} + \lVert \omega_{3} \rVert_{L^{p}}^{r} d\tau < \infty, \hspace{3mm} \frac{3}{p} + \frac{2}{r} \leq 1 + \frac{1}{p}, \hspace{3mm} 2 \leq p \leq \infty
\end{equation}
for the MHD system in [15] (cf. Corollary 1.3 [25]). In fact, (9) is an immediate consequence of obtaining a condition in terms of $\partial_{1}u_{1}, \partial_{2}u_{2}, \partial_{3}u_{3}$ (cf. also [20]) and using incompressibility as well as the continuity of Riesz transform that shows that (e.g. [22] Lemma 2.1)

\begin{equation}
\lVert \partial_{i} u_{j} \rVert_{L^{p}} \lesssim_{p} ( \lVert \partial_{3} u_{3} \rVert_{L^{p}} + \lVert \omega_{3} \rVert_{L^{p}}), \hspace{3mm} 1 < p < \infty, i, j = 1, 2.
\end{equation}

Let us motivate the study of this manuscript specifically. We are interested in the regularity criterion for the MHD system involving only $u_{3}$ and $\omega_{3}$, which is not accessible as long as one relies on an inequality of the type such as (10). In fact, the authors in [21] obtained a condition for the NSE in terms of $u_{3}, \partial_{1}u_{2}, \partial_{2}u_{1}$, all in scaling invariant norms, and suggested an open problem to replace this condition with $u_{3}, \omega_{3}$ (see Theorem 1 and Remark 1 [21]). Subsequently, the authors in [12] obtained a condition for the NSE involving only $\partial_{3}u_{3}, \omega_{3}$ both in scaling invariant norms but remarked that replacing $\partial_{3}u_{3}$ by $u_{3}$ seems difficult (see Theorem 1.1 and Remark 1.2 [12]). Let us now present our result: 

\begin{theorem}

Suppose $(u_{0}, b_{0}) \in L^{2}(\mathbb{R}^{3}), \nabla\cdot u_{0} = \nabla\cdot b_{0} = 0$ and $(u,b)$ is a weak solution pair to the MHD system (1a)-(1c) on $[0,T)$. If 

\begin{equation}
\int_{0}^{T} \lVert u_{3} \rVert_{L^{p_{1}}}^{r_{1}} + \lVert \omega_{3} \rVert_{L^{p_{2}}}^{r_{2}} d\tau < \infty, \hspace{3mm} 
\begin{cases}
\frac{3}{p_{1}} + \frac{2}{r_{1}} \leq \frac{4}{9} - \frac{1}{3p_{1}}, & \frac{15}{2} \leq p_{1} <  \infty, \\
\frac{3}{p_{2}} + \frac{2}{r_{2}} \leq 2,  & \frac{3}{2} < p_{2} < \infty,
\end{cases}
\end{equation}
then the solution is smooth on $(0,T)$. Moreover, the condition on $\omega_{3}$ may be replaced by $\sup_{t \in [0,T)} \lVert \omega_{3}(t) \rVert_{L^{\frac{3}{2}}}$ being sufficiently small. 

\end{theorem}

An immediate corollary also of much interest is the following: 

\begin{corollary}
	
Suppose $(u_{0}, b_{0}) \in L^{2}(\mathbb{R}^{3}), \nabla\cdot u_{0} = \nabla\cdot b_{0} = 0$ and $(u,b)$ is a weak solution pair to the MHD system (1a)-(1c) on $[0,T)$. If 
	
\begin{equation*}
\int_{0}^{T} \lVert u_{3} \rVert_{L^{p_{1}}}^{r_{1}} + \lVert \partial_{1}u_{2} \rVert_{L^{p_{2}}}^{r_{2}}  + \lVert \partial_{2}u_{1} \rVert_{L^{p_{2}}}^{r_{2}}d\tau < \infty, \hspace{3mm} 
\begin{cases}
\frac{3}{p_{1}} + \frac{2}{r_{1}} \leq \frac{4}{9} - \frac{1}{3p_{1}}, & \frac{15}{2} \leq p_{1} <  \infty, \\
\frac{3}{p_{2}} + \frac{2}{r_{2}} \leq 2,  & \frac{3}{2} < p_{2} < \infty,
\end{cases}
\end{equation*}
then the solution is smooth on $(0,T)$. Moreover, the condition on $\partial_{1}u_{2}, \partial_{2}u_{1}$ may be replaced by $\sup_{t \in [0,T)} \lVert \partial_{1}u_{2}(t) \rVert_{L^{\frac{3}{2}}} + \lVert \partial_{2}u_{1}(t) \rVert_{L^{\frac{3}{2}}}$ being sufficiently small. 
	
\end{corollary}

\begin{remark}
\begin{enumerate}
\item The condition on $p_{2}, r_{2}$ allows the scaling invariant level. This is due to the new decomposition of the four non-linear terms in Proposition 1.3.

\item Because the problem raised by the authors in [21], namely a condition in terms of $u_{3}, \omega_{3}$ both in scaling invariant norms for the NSE, is open, it is expected that such a result for the MHD system is much more difficult. With this in mind, the only difference with the desired result and that of Theorem 1.1 is that the condition on $u_{3}$ is not at the scaling invariant level. To the best of the author's knowledge, this is the first regularity criterion for the MHD system involving only $u_{3}, \omega_{3}$. 

\item As we discussed, a regularity criterion for the MHD system in terms of only $u_{3}$ is not known. Theorem 1.1 is precisely that, only added by $\omega_{3}$ at a scaling invariant level. Similarly, a regularity criterion for the MHD system (1a)-(1c) in terms of only $\omega_{2}, \omega_{3}$ is not known. Theorem 1.1 is precisely that, with $\omega_{2}$ replaced by $u_{3}$; furthermore, $\omega_{3}$ is allowed to be at the scaling invariant level. 

\item Concerning Corollary 1.2, e.g. in [15] the authors obtained a regularity criterion in terms of $\partial_{i}u_{1}, \partial_{j}u_{2}, \partial_{k}u_{3}$ for any $i, j, k \in \{1,2,3\}$; however, none of them were in scaling invariant spaces. While Corollary 1.3 has $u_{3}$ instead its partial derivative, remarkably the other two partial derivatives, $\partial_{1}u_{2}, \partial_{2}u_{1}$ are allowed to be in scaling invariant norms. 

\end{enumerate}
\end{remark}

We now elaborate on the proof of Theorem 1.1. In [29], the author obtained a regularity criterion of the MHD system in terms of only $u_{3}, j_{3}$ with the latter in a scaling invariant norm. On the other hand, in [26] the author initiated a series of estimates that essentially \emph{controls} $b_{i}$ in terms of $u_{i}, \forall \hspace{1mm} i \in \{1,2,3\}$ as follows:  

\begin{equation}
\sup_{\tau\in [0,t]} \lVert b_{i}(\tau)\rVert_{L^{p}}^{2} \leq \lVert b_{i}(0)\rVert_{L^{p}}^{2} +  c(p)\int_{0}^{t} \lVert \nabla b(\lambda)\rVert_{L^{2}}^{2} \lVert u_{i}(\lambda)\rVert_{L^{\frac{6p}{6-p}}}^{2} d\lambda
\end{equation}
for any $t \in [0,T]$ where $\frac{6p}{6-p} = \infty$ if $p = 6$ (cf. [27, 28]). This is essentially due to the fact that the $i-$th component of $b$ has relatively a simple form of 

\begin{equation*}
\partial_{t} b_{i} + (u\cdot\nabla) b_{i} = \eta \Delta b_{i} + (b\cdot\nabla) u_{i}
\end{equation*}
so that upon the $L^{p}$-estimate, the first non-linear term vanishes while in the second, $u_{i}$ is already \emph{separated}. Given this idea and the result in [29], one is tempted to hope that perhaps $j_{3}$ may be \emph{controlled} by $\omega_{3}$; this seems very difficult as $j_{3}$ is governed by the equation of 

\begin{equation*}
\partial_{t} j_{3} + (u\cdot\nabla) j_{3} - (b\cdot\nabla) \omega_{3} - \Delta j_{3} = (j\cdot\nabla) u_{3} - (\omega\cdot\nabla) b_{3} + 2[\partial_{1} b \cdot\partial_{2} u - \partial_{2} b \cdot\partial_{1} u].
\end{equation*}
Only term that vanishes upon an $L^{p}$-estimate is the first non-linear term and the last term in the bracket in particular seems very difficult to even find $\omega_{3}$. We also remark that even if an analogue of (12) for $j_{3}$ may be obtained, applying it within the proof of [29] implies that the criterion on $\omega_{3}$ will not be in a scaling invariant norm. The novelty of this manuscript is the following new  decomposition which requires many cancellations observed in [29] as well as appropriate application of new identities in (16):

\begin{proposition}
Let smooth solutions of the MHD system (1a)-(1c),
	
\begin{equation}
\begin{split}
&\int (u\cdot\nabla) u \cdot \Delta_{h} u - (b\cdot\nabla) b \cdot \Delta_{h} u + (u\cdot\nabla) b \cdot \Delta_{h} b - (b\cdot\nabla) u \cdot \Delta_{h} b\\
\lesssim& \int \lvert u_{3} \rvert ( \lvert \nabla u \rvert \lvert \nabla\nabla_{h} u\rvert + \lvert \nabla b \rvert \lvert \nabla\nabla_{h} b\rvert) + \lvert b_{3} \rvert ( \lvert \nabla u\rvert \lvert \nabla\nabla_{h} b\rvert + \lvert \nabla b \rvert \lvert \nabla\nabla_{h} u\rvert) \\
&+ \lvert \nabla_{h} b \rvert^{2} \lvert \nabla_{h}^{2} \Delta_{h}^{-1} \omega_{3} \rvert + \lvert   \nabla_{h}^{2} \Delta_{h}^{-1} u_{3} \rvert \lvert \nabla\nabla_{h} b \rvert \lvert \nabla_{h} b\rvert. 
\end{split}
\end{equation}
	
\end{proposition} 

As it will be clear in (29), it is crucial that we have $\nabla_{h}^{2}\Delta_{h}^{-1}$ in (13) instead of e.g. $\partial_{3}\partial_{1}\Delta_{h}^{-1}$ (see equations (22), (23)). We also remark that the identities (16) can be proven via Fourier analysis. With this in mind, (13) is an interesting decomposition that combined the use of both integration by parts as done in (6) by the authors in [5, 17] and identities from Fourier analysis.  

In the next section, let us set up further notations, state a few useful identities and inequalities. Thereafter, we prove Theorem 1.1. 

\section{Preliminaries}

Let us denote for brevity 

\begin{equation}
\begin{split}
W(t) \triangleq  (\lVert \nabla_{h} u\rVert_{L^{2}}^{2} + \lVert \nabla_{h} b\rVert_{L^{2}}^{2})(t), \hspace{2mm} & X(t) \triangleq  (\lVert \nabla u\rVert_{L^{2}}^{2}+ \lVert \nabla b\rVert_{L^{2}}^{2})(t),\\
Y(t) \triangleq (\lVert \nabla\nabla_{h} u\rVert_{L^{2}}^{2} + \lVert \nabla\nabla_{h} b\rVert_{L^{2}}^{2})(t), \hspace{2mm}  & Z(t) \triangleq (\lVert \Delta u\rVert_{L^{2}}^{2} + \lVert \Delta b\rVert_{L^{2}}^{2})(t).  
\end{split}
\end{equation}
The following is a special case of Troisi's inequality (cf. [6] for proof):

\begin{equation}
\lVert f\rVert_{L^{6}} \lesssim \lVert \partial_{1} f\rVert_{L^{2}}^{\frac{1}{3}} \lVert \partial_{2} f\rVert_{L^{2}}^{\frac{1}{3}} \lVert \partial_{3} f\rVert_{L^{2}}^{\frac{1}{3}}. 
\end{equation}
The following identity was utilized for the NSE in [9] and thereafter for the MHD system in [30]: $\forall f = (f_{h}, f_{3})$ such that $\nabla \cdot f = 0$, 
	
\begin{equation}
\begin{split}
f_{1} = -\partial_{2} \Delta_{h}^{-1} (\nabla \times f)\cdot e_{3} - \partial_{1} \Delta_{h}^{-1} \partial_{3}f_{3},\\
f_{2} = \partial_{1} \Delta_{h}^{-1} (\nabla \times f)\cdot e_{3} - \partial_{2} \Delta_{h}^{-1} \partial_{3} f_{3}.
\end{split}
\end{equation}
The following lemma was initiated in [26, 27, 28, 29] and this is directly due to [31]: 

\begin{lemma}
For smooth solutions of the MHD system (1a)-(1c), for any $i \in \{1,2,3\}, 2 \leq p \leq 6, 0 < t_{1} < t_{2}$, 

\begin{equation*}
\sup_{t_{1} \leq t \leq t_{2}} \lVert b_{i}(t) \rVert_{L^{p}}^{2}  \leq \lVert b_{i}(t_{1}) \rVert_{L^{p}}^{2} + c \int_{t_{1}}^{t_{2}} \lVert \nabla_{h} b\rVert_{L^{2}}^{\frac{4}{3}} \lVert \nabla b\rVert_{L^{2}}^{\frac{2}{3}} \lVert u_{i} \rVert_{L^{\frac{6p}{6-p}}}^{2} d\tau 
\end{equation*}
where $\frac{6p}{6-p} = \infty$ if $p = 6$. 

\end{lemma}
	
An inequality similar to the following was used in many places; its formal proof may be found in [31]:  

\begin{lemma}
For smooth solutions of the MHD system (1a)-(1c), 
\begin{equation*}
\begin{split}
& \int (u\cdot\nabla) u \cdot \Delta u - (b\cdot\nabla) b \cdot \Delta u + (u\cdot\nabla) b \cdot \Delta b - (b\cdot\nabla) u \cdot \Delta b\\
\lesssim&  \int ( \lvert \nabla_{h} u\rvert + \lvert \nabla_{h} b\rvert) ( \lvert \nabla u\rvert^{2} + \lvert \nabla b\rvert^{2}).
\end{split}
\end{equation*}
\end{lemma}	 	
	
\section{Proof of Proposition 1.3}

In this section we prove Proposition 1.3. Although we make use of many cancellations described with detail in [29], for completeness we sketch the steps. Firstly, we integrate by parts and use divergence-free conditions (1c) to obtain 

\begin{equation*}
\begin{split}
&\int (u\cdot\nabla) u \cdot \Delta_{h} u - (b\cdot\nabla) b \cdot \Delta_{h} u + (u\cdot\nabla) b \cdot \Delta_{h} b - (b\cdot\nabla) u \cdot \Delta_{h} b\\
=& \sum_{i,j=1}^{3} \sum_{k=1}^{2} \int -\partial_{k} u_{i} \partial_{i} u_{j} \partial_{k} u_{j} + \partial_{k} b_{i} \partial_{i} b_{j} \partial_{k} u_{j} - \partial_{k} u_{i} \partial_{i} b_{j} \partial_{k} b_{j} + \partial_{k} b_{i} \partial_{i} u_{j} \partial_{k} b_{j} \\
\triangleq& I + II + III + IV.
\end{split}
\end{equation*}

It is shown in the equations (3.5)-(3.10) of [29] (cf. [17] Lemma 2.3) that 

\begin{equation}
I \lesssim \int \lvert u_{3} \rvert \lvert \nabla u\rvert \lvert \nabla\nabla_{h} u\rvert. 
\end{equation}

Moreover, it is shown in equations (3.11)-(3.13) of [29] that 

\begin{equation}
\begin{split}
II =& \int (\partial_{1}b_{1})^{2} \partial_{1} u_{1} + \partial_{2} b_{1} \partial_{1} b_{1} \partial_{2} u_{1} + \partial_{1} b_{1} \partial_{1} b_{2} \partial_{1} u_{2} + \partial_{2} b_{1} \partial_{1} b_{2} \partial_{2} u_{2}\\
&+ \partial_{1} b_{2} \partial_{2} b_{1} \partial_{1}u_{1} + \partial_{2} b_{2}\partial_{2} b_{1} \partial_{2} u_{1} + \partial_{1} b_{2} \partial_{2} b_{2} \partial_{1} u_{2} + (\partial_{2} b_{2})^{2} \partial_{2} u_{2} \\
&+ \sum_{j,k=1}^{2} \int \partial_{k} b_{3} \partial_{3} b_{j} \partial_{k} u_{j} + \sum_{i=1}^{3} \sum_{k=1}^{2} \int \partial_{k} b_{i} \partial_{i} b_{3} \partial_{k} u_{3} \\
\leq& \sum_{i=1}^{8} II_{i} + c \int \lvert b_{3} \rvert ( \lvert \nabla\nabla_{h} b \rvert \lvert \nabla u \rvert + \lvert \nabla b \rvert \lvert \nabla\nabla_{h} u\rvert) + \lvert u_{3} \rvert \lvert \nabla b\rvert \lvert \nabla\nabla_{h} b \rvert,
\end{split}
\end{equation} 

\begin{equation}
\begin{split}
III=& -\int \partial_{1} u_{1} (\partial_{1} b_{1})^{2} + \partial_{2} u_{1} \partial_{1} b_{1} \partial_{2} b_{1} + \partial_{1} u_{1} (\partial_{1} b_{2})^{2} + \partial_{2} u_{1} \partial_{1} b_{2} \partial_{2} b_{2}\\
&+ \partial_{1} u_{2} \partial_{2} b_{1} \partial_{1} b_{1} + \partial_{2} u_{2} \partial_{2} b_{1} \partial_{2} b_{1} + \partial_{1} u_{2} \partial_{2} b_{2} \partial_{1} b_{2} + \partial_{2} u_{2} (\partial_{2} b_{2})^{2} \\
&- \sum_{j,k=1}^{2} \int \partial_{k} u_{3} \partial_{3} b_{j} \partial_{k} b_{j} - \sum_{i=1}^{3} \sum_{k=1}^{2} \int \partial_{k} u_{i} \partial_{i} b_{3} \partial_{k} b_{3}\\
\leq& \sum_{i=1}^{8} III_{i} + c \int \lvert u_{3} \rvert \lvert \nabla\nabla_{h} b \rvert \lvert \nabla b\rvert + \lvert b_{3} \rvert ( \lvert \nabla\nabla_{h} u\rvert \lvert \nabla b\rvert + \lvert \nabla u\rvert \lvert \nabla\nabla_{h} b\rvert), 
\end{split}
\end{equation}

\begin{equation}
\begin{split}
IV =& \int (\partial_{1} b_{1})^{2} \partial_{1} u_{1} + \partial_{2} b_{1} \partial_{1} u_{1} \partial_{2} b_{1} + \partial_{1} b_{1} \partial_{1} u_{2} \partial_{1} b_{2} + \partial_{2} b_{1} \partial_{1} u_{2} \partial_{2} b_{2}\\
&+ \partial_{1} b_{2} \partial_{2} u_{1} \partial_{1} b_{1} + \partial_{2} b_{2} \partial_{2} u_{1} \partial_{2} b_{1} + \partial_{1} b_{2} \partial_{2} u_{2} \partial_{1} b_{2} + (\partial_{2} b_{2})^{2} \partial_{2} u_{2} \\
&+ \sum_{j,k=1}^{2} \int \partial_{k} b_{3} \partial_{3} u_{j} \partial_{k} b_{j} + \sum_{i=1}^{3} \sum_{k=1}^{2} \int \partial_{k} b_{i} \partial_{i} u_{3} \partial_{k} b_{3} \\
\leq& \sum_{i=1}^{8} IV_{i} + c \int \lvert b_{3} \rvert ( \lvert \nabla\nabla_{h} u\rvert \lvert \nabla b\rvert + \lvert \nabla u\rvert \lvert \nabla\nabla_{h} b\rvert). 
\end{split}
\end{equation}

Hence, we have from (17)-(20) 

\begin{equation}
\begin{split}
& \int (u\cdot\nabla) u \cdot \Delta_{h} u - (b\cdot\nabla) b \cdot \Delta_{h} u + (u\cdot\nabla) b \cdot \Delta_{h} b - (b\cdot\nabla) u \cdot \Delta_{h} b\\
\leq& c \int \lvert u_{3} \rvert ( \lvert \nabla u\rvert \lvert \nabla \nabla_{h} u\rvert + \lvert \nabla b\rvert \lvert \nabla\nabla_{h} b\rvert) + \lvert b_{3} \rvert ( \lvert \nabla u\rvert \lvert \nabla\nabla_{h} b\rvert + \lvert \nabla b\rvert \lvert \nabla\nabla_{h} u\rvert)\\
&+ \sum_{i=1}^{8} II_{i} + III_{i} + IV_{i}. 
\end{split}
\end{equation}

The following cancellations have been discovered in [29]:

\begin{equation*}
\begin{split}
II_{1} + III_{1} =& II_{2} + III_{2} = II_{7} + III_{7} = II_{8} + III_{8} = 0,\\
II_{3} + IV_{3} =&  \int 2 \partial_{1}b_{1}\partial_{1}b_{2} \partial_{1}u_{2}, \hspace{1mm} II_{4} + II_{5} 
=  \int u_{3}\partial_{3}(\partial_{2}b_{1}\partial_{1}b_{2}),\\
II_{6} + IV_{6} =& \int 2 \partial_{2}b_{2}\partial_{2}b_{1}\partial_{2}u_{1},\\
III_{3} + IV_{7} =& -\int u_{3}\partial_{3}(\partial_{1}b_{2})^{2} + \int 2\partial_{2}u_{2}(\partial_{1}b_{2})^{2},\nonumber\\
III_{4} + IV_{5} =& -\int b_{3}\partial_{3}(\partial_{2}u_{1}\partial_{1}b_{2}) + \int 2\partial_{2}u_{1}\partial_{1}b_{2}\partial_{1}b_{1},\nonumber\\
III_{5} + IV_{4} =&  - \int b_{3}\partial_{3}(\partial_{1}u_{2}\partial_{2}b_{1}) + \int 2\partial_{1}u_{2}\partial_{2}b_{1}\partial_{2}b_{2},\nonumber\\
III_{6} + IV_{2} =& -\int u_{3}\partial_{3}(\partial_{2}b_{1})^{2} + \int 2\partial_{1}u_{1}(\partial_{2}b_{1})^{2},\nonumber\\
IV_{1} +IV_{8} =& \int b_{3}\partial_{3}(\partial_{1}b_{1}\partial_{1}u_{1} + \partial_{2}b_{2}\partial_{2}u_{2}) - \int u_{3}\partial_{3}(\partial_{1}b_{1}\partial_{2}b_{2}),\nonumber
\end{split}
\end{equation*}
and therefore

\begin{equation}
\begin{split}
&\sum_{i=1}^{8} II_{i} + III_{i} + IV_{i}\\
\leq& c \int \lvert u_{3} \rvert \lvert \nabla\nabla_{h} b \rvert \lvert \nabla b\rvert + \lvert b_{3} \rvert ( \lvert \nabla\nabla_{h} u\rvert \lvert \nabla b\rvert + \lvert \nabla u\rvert \lvert \nabla\nabla_{h} b\rvert) \\
&+ 2 \int \partial_{1} b_{1} \partial_{1} b_{2} \partial_{1} u_{2} + \partial_{2} b_{2} \partial_{2} b_{1} \partial_{2} u_{1} + \partial_{2} u_{2} (\partial_{1} b_{2})^{2}\\
& \hspace{6mm} + \partial_{2} u_{1} \partial_{1} b_{2} \partial_{1} b_{1} + \partial_{1} u_{2} \partial_{2} b_{1} \partial_{2} b_{2} + \partial_{1} u_{1} (\partial_{2} b_{1})^{2}. 
\end{split}
\end{equation}
(see equations (3.14)-(3.25) of [29]). We apply the identity (16) on $u_{1}, u_{2}$ and integrate by parts so that the second integral in (22) can be written as  

\begin{equation}
\begin{split}
& \int \partial_{1} b_{1} \partial_{1} b_{2} \partial_{1} u_{2} + \partial_{2} b_{2} \partial_{2} b_{1} \partial_{2} u_{1} + \partial_{2} u_{2} (\partial_{1} b_{2})^{2}\\
&+ \partial_{2} u_{1} \partial_{1} b_{2} \partial_{1} b_{1} + \partial_{1} u_{2} \partial_{2} b_{1} \partial_{2} b_{2} + \partial_{1} u_{1} (\partial_{2} b_{1})^{2} \\
=& \int \partial_{1} b_{1} \partial_{1} b_{2} \partial_{11}^{2} \Delta_{h}^{-1} \omega_{3} + \int \partial_{3} (\partial_{1} b_{1} \partial_{1} b_{2}) \partial_{12} \Delta_{h}^{-1} u_{3}\\
&- \int \partial_{2} b_{2} \partial_{2} b_{1} \partial_{22}^{2} \Delta_{h}^{-1} \omega_{3} + \int \partial_{3} (\partial_{2} b_{2} \partial_{2} b_{1}) \partial_{21} \Delta_{h}^{-1} u_{3} \\
&+ \int \partial_{21} \Delta_{h}^{-1} \omega_{3} (\partial_{1} b_{2})^{2} + \int \partial_{22}^{2} \Delta_{h}^{-1} u_{3} \partial_{3} (\partial_{1} b_{2})^{2}\\
&- \int \partial_{22}^{2} \Delta_{h}^{-1} \omega_{3} \partial_{1} b_{2} \partial_{1} b_{1} + \int \partial_{21} \Delta_{h}^{-1} u_{3} \partial_{3} (\partial_{1} b_{2} \partial_{1} b_{1})\\
&+ \int \partial_{11}^{2} \Delta_{h}^{-1} \omega_{3} \partial_{2} b_{1} \partial_{2} b_{2} + \int \partial_{12} \Delta_{h}^{-1} u_{3} \partial_{3} (\partial_{2} b_{1} \partial_{2} b_{2})\\
&- \int \partial_{12} \Delta_{h}^{-1} \omega_{3} (\partial_{2} b_{1})^{2} + \int \partial_{11}^{2} \Delta_{h}^{-1} u_{3} \partial_{3}(\partial_{2} b_{1})^{2}. 
\end{split}
\end{equation}
Applying (23) in (22), we obtain 

\begin{equation}
\begin{split}
\sum_{i=1}^{8} II_{i} + III_{i} + IV_{i}
\lesssim&  \int \lvert u_{3} \rvert \lvert \nabla\nabla_{h} b \rvert \lvert \nabla b\rvert + \lvert b_{3} \rvert ( \lvert \nabla\nabla_{h} u\rvert \lvert \nabla b\rvert + \lvert \nabla u\rvert \lvert \nabla\nabla_{h} b\rvert) \\
&+  \lvert \nabla_{h} b\rvert^{2} \lvert \nabla_{h}^{2} \Delta_{h}^{-1} \omega_{3} \rvert + \lvert \nabla_{h}^{2} \Delta_{h}^{-1} u_{3} \rvert \lvert \nabla\nabla_{h} b \rvert \lvert \nabla_{h} b\rvert. 
\end{split}
\end{equation}
Thus, (24) along with (21) completes the proof of Proposition 1.3.  

\section{Proof of Theorem 1.1}

We follow the method in [31]. We first fix $\delta \in (0,T)$ arbitrarily. We know $\exists$ at least one weak solution pair $(u,b) \in L^{\infty} ((0, T) ; L^{2}(\mathbb{R}^{3})) \cap L^{2}(0, T; \dot{H}^{1}(\mathbb{R}^{3}))$. Due to the local existence of the unique strong solution, restarting at time $\hat{\delta} \in (0, \delta), \exists \hspace{1mm} ! \hspace{1mm} 
\tilde{u}, \tilde{b} \in C([\hat{\delta}, T^{\ast} ); \dot{H}^{1}(\mathbb{R}^{3})) \cap L^{2} ( [\hat{\delta}, T^{\ast} ); \dot{H}^{2}(\mathbb{R}^{3}))$ where $[\hat{\delta}, T^{\ast})$ is the life span of the unique strong solution; moreover, it is well known that this regularity leads to $\tilde{u}, \tilde{b} \in C^{\infty} (\mathbb{R}^{3} \times (\hat{\delta}, T^{\ast}))$. Because the strong solution is the only weak solution, $u = \tilde{u}, b = \tilde{b} \text{ on } [\hat{\delta}, T^{\ast})$. If $T^{\ast} \geq T$, then we obtain $u, b \in C^{\infty} (\mathbb{R}^{3} \times (0, T))$. Suppose $T^{\ast} < T$; thus, necessarily $\limsup_{t \to T^{\ast}} (\lVert \nabla u\rVert_{L^{2}}^{2} + \lVert \nabla b\rVert_{L^{2}}^{2})(t) = \infty$. We show that $\forall t < T^{\ast}$, $(\lVert \nabla u\rVert_{L^{2}}^{2} + \lVert \nabla b\rVert_{L^{2}}^{2})(t) \leq c$, a contradiction to the definition of $T^{\ast}$. 

We choose $\tilde{\epsilon} > 0$ to be precisely determined subsequently and then select $\Gamma < T^{\ast}$ sufficiently close to $T^{\ast}$ so that 

\begin{equation}
\forall \hspace{1mm} t \in [\Gamma, T^{\ast}) \hspace{2mm} \int_{\Gamma}^{t} X d\tau < \tilde{\epsilon}. 
\end{equation} 

Let us consider the case $p_{1} \in (\frac{15}{2}, \infty)$ for simplicity of presentation as the case $p_{1} = \frac{15}{2}$ requires only a straight-forward modification. We also fix $p_{2} \in (\frac{3}{2}, \infty)$ only to make a remark about how the case $p_{2} = \frac{3}{2}$ may be obtained afterwards. Now we define $p_{0} \triangleq \frac{6p_{1}}{6+ p_{1}}$ so that $p_{0} \in (\frac{10}{3}, 6)$ and take $L^{2}$-inner products on (1a)-(1b) with $(-\Delta_{h} u, -\Delta_{h} b)$ to obtain 

\begin{equation}
\begin{split}
& \frac{1}{2} \partial_{t} W + Y \\
\lesssim& \int \lvert u_{3} \rvert ( \lvert \nabla u \rvert \lvert \nabla\nabla_{h} u\rvert + \lvert \nabla b \rvert \lvert \nabla\nabla_{h} b\rvert) + \int \lvert b_{3} \rvert ( \lvert \nabla u\rvert \lvert \nabla\nabla_{h} b\rvert + \lvert \nabla b \rvert \lvert \nabla\nabla_{h} u\rvert) \\
&+  \int \lvert \nabla_{h} b\rvert^{2} \lvert \nabla_{h}^{2} \Delta_{h}^{-1} \omega_{3} \rvert + \lvert \nabla_{h}^{2} \Delta_{h}^{-1} u_{3} \rvert \lvert \nabla\nabla_{h} b \rvert \lvert \nabla_{h} b\rvert  \triangleq V + VI + VII 
\end{split}
\end{equation}
due to Proposition 1.3. For $\epsilon > 0$ arbitrary small we may estimate $V$ and $VI$ by 

\begin{equation}
\begin{split}
&V + VI \\
\lesssim& \lVert u_{3} \rVert_{L^{p_{1}}}( \lVert \nabla u\rVert_{L^{\frac{2p_{1}}{p_{1} - 2}}}\lVert \nabla\nabla_{h} u\rVert_{L^{2}} + \lVert \nabla b\rVert_{L^{\frac{2p_{1}}{p_{1} - 2}}}\lVert \nabla\nabla_{h} b\rVert_{L^{2}}) \\
&+ \lVert b_{3} \rVert_{L^{p_{0}}} (\lVert \nabla u\rVert_{L^{\frac{2p_{0}}{p_{0} - 2}}}\lVert \nabla\nabla_{h} b\rVert_{L^{2}} +  \lVert \nabla b\rVert_{L^{\frac{2p_{0}}{p_{0} - 2}}}\lVert \nabla\nabla_{h} u\rVert_{L^{2}})\\
\lesssim& \lVert u_{3} \rVert_{L^{p_{1}}}( \lVert \nabla u\rVert_{L^{2}}^{\frac{p_{1} - 3}{p_{1}}}\lVert \nabla u\rVert_{L^{6}}^{\frac{3}{p_{1}}}\lVert \nabla\nabla_{h} u\rVert_{L^{2}} + \lVert \nabla b\rVert_{L^{ 2}}^{\frac{p_{1} - 3}{p_{1}}}\lVert \nabla b\rVert_{L^{6}}^{\frac{3}{p_{1}}}\lVert \nabla\nabla_{h} b\rVert_{L^{2}}) \\
&+ \lVert b_{3} \rVert_{L^{p_{0}}} (\lVert \nabla u\rVert_{L^{2}}^{\frac{p_{0} - 3}{p_{0}}} \lVert \nabla u\rVert_{L^{6}}^{\frac{3}{p_{0}}}\lVert \nabla\nabla_{h} b\rVert_{L^{2}} +  \lVert \nabla b\rVert_{L^{ 2}}^{\frac{p_{0} - 3}{p_{0}}}\lVert \nabla b\rVert_{L^{6}}^{\frac{3}{p_{0}}}\lVert \nabla\nabla_{h} u\rVert_{L^{2}})\\
\leq&  \epsilon Y + c\left( 
\lVert u_{3} \rVert_{L^{p_{1}}}^{\frac{2p_{1}}{p_{1} - 2}} X^{\frac{p_{1} - 3}{p_{1} - 2}}Z^{\frac{1}{p_{1} - 2}}  + \lVert b_{3} \rVert_{L^{p_{0}}}^{\frac{2p_{0}}{p_{0} - 2}} X^{\frac{p_{0} - 3}{p_{0} - 2}} Z^{\frac{1}{p_{0} - 2}}\right) 
\end{split}
\end{equation}
due to H$\ddot{o}$lder's and interpolation inequalities, (15), and Young's inequalities. We now consider $VII$ and estimate 

\begin{equation}
\begin{split}
VII \lesssim \lVert \nabla_{h} b\rVert_{L^{\frac{2p_{2}}{p_{2} - 1}}}^{2} \lVert \omega_{3} \rVert_{L^{p_{2}}} + \lVert u_{3} \rVert_{L^{p_{1}}} \lVert \nabla\nabla_{h} b\rVert_{L^{2}} \lVert \nabla_{h} b\rVert_{L^{\frac{2p_{1}}{p_{1} - 2}}}
\end{split}
\end{equation}
by H$\ddot{o}$lder's inequalities and the continuity of Riesz transform in the following anisotropic way: $\forall \hspace{1mm} p \in (1,\infty)$ 

\begin{equation}
\lVert \nabla_{h}^{2} \Delta_{h}^{-1} f \rVert_{L^{p}} = \lVert \lVert \nabla_{h}^{2} \Delta_{h}^{-1} f\rVert_{L_{h}^{p}} \rVert_{L_{v}^{p}} \lesssim_{p} \lVert \lVert f\rVert_{L_{h}^{p}} \rVert_{L_{v}^{p}} \approx_{p} \lVert f\rVert_{L^{p}}.
\end{equation} 

We further estimate from (28)

\begin{equation}
\begin{split}
VII 
\lesssim& \lVert \nabla_{h} b\rVert_{L^{2}}^{\frac{2p_{2} - 3}{p_{2}}}\lVert \nabla\nabla_{h} b\rVert_{L^{2}}^{\frac{3}{p_{2}}}\lVert \omega_{3} \rVert_{L^{p_{2}}}\\
&+ \lVert u_{3} \rVert_{L^{p_{1}}} \lVert \nabla\nabla_{h} b\rVert_{L^{2}}^{\frac{p_{1} + 2}{p_{1}}}\lVert \nabla_{h} b\rVert_{L^{2}}^{\frac{p_{1} - 3}{p_{1}}}\lVert \Delta b\rVert_{L^{2}}^{\frac{1}{p_{1}}}\\
\leq& \epsilon Y + c \left(\lVert \omega_{3} \rVert_{L^{p_{2}}}^{\frac{2p_{2}}{2p_{2} - 3}}W + \lVert u_{3} \rVert_{L^{p_{1}}}^{\frac{2p_{1}}{p_{1} - 2}}W^{\frac{p_{1} - 3}{p_{1} - 2}}Z^{\frac{1}{p_{1} - 2}} \right) 
\end{split}
\end{equation}
by Gagliardo-Nirenberg and interpolation inequalities, (15) and Young's inequalities. 

Applying (27) and (30) in (26), absorbing $2\epsilon Y$ for $\epsilon > 0$ sufficiently small, Gronwall's type argument using 

\begin{equation*}
1 \leq \sup_{\lambda \in [\Gamma, \tau]}e^{c \int_{\lambda}^{\tau} \lVert \omega_{3} \rVert_{L^{p_{2}}}^{\frac{2p_{2}}{2p_{2} - 3}} d\phi } \lesssim e^{c \int_{0}^{T^{\ast}} \lVert \omega_{3} \rVert_{L^{p_{2}}}^{\frac{2p_{2}}{2p_{2} - 3}} d\phi } \lesssim 1
\end{equation*}
due to (11) leads to, for every $\tau \in [\Gamma, t]$

\begin{equation}
\begin{split}
W(\tau) + \int_{\Gamma}^{\tau} Y d\lambda 
\lesssim 1 + \int_{\Gamma}^{t} \lVert u_{3} \rVert_{L^{p_{1}}}^{\frac{2p_{1}}{p_{1} - 2}} X^{\frac{p_{1} - 3}{p_{1} - 2}}Z^{\frac{1}{p_{1} - 2}}  + \lVert b_{3} \rVert_{L^{p_{0}}}^{\frac{2p_{0}}{p_{0} - 2}} X^{\frac{p_{0} - 3}{p_{0} - 2}} Z^{\frac{1}{p_{0} - 2}}  d\lambda.
\end{split}
\end{equation}
We take $\sup_{\tau \in [\Gamma, t]}$ on the left hand side and continue this bound as follows: 

\begin{equation}
\begin{split}
&\sup_{\tau \in [\Gamma, t]}W (\tau) + \int_{\Gamma}^{t}Y d\tau \\
\lesssim& 1 + \int_{\Gamma}^{t}\lVert u_{3} \rVert_{L^{p_{1}}}^{\frac{2p_{1}}{p_{1} - 2}} X^{\frac{p_{1} - 3}{p_{1} - 2}}Z^{\frac{1}{p_{1} - 2}}d\tau+ \sup_{\tau \in [\Gamma, t]} \lVert b_{3}(\tau) \rVert_{L^{p_{0}}}^{\frac{2p_{0}}{p_{0} - 2}} \int_{\Gamma}^{t}X^{\frac{p_{0} - 3}{p_{0} - 2}}Z^{\frac{1}{p_{0} - 2}}d\tau\\
\lesssim& 1 + \int_{\Gamma}^{t}\lVert u_{3} \rVert_{L^{p_{1}}}^{\frac{2p_{1}}{p_{1} - 2}} X^{\frac{p_{1} - 3}{p_{1} - 2}}Z^{\frac{1}{p_{1} - 2}} d\tau\\
+& \left(\int_{\Gamma}^{t} \lVert \nabla_{h} b\rVert_{L^{2}}^{\frac{4}{3}} \lVert \nabla b\rVert_{L^{2}}^{\frac{2}{3}} \lVert u_{3} \rVert_{L^{\frac{6p_{0}}{6-p_{0}}}}^{2} d\tau   \right)^{\frac{p_{0}}{p_{0} - 2}} \int_{\Gamma}^{t}X^{\frac{p_{0} - 3}{p_{0} - 2}}Z^{\frac{1}{p_{0} - 2}} d\tau
\end{split}
\end{equation}
by Lemma 2.1 as $p_{0} \in (\frac{10}{3},6)$, an elementary inequality of $(a+b)^{p} \leq 2^{p}(a^{p} + b^{p}), \text{ for } 0 \leq p < \infty \text{ and } a, b \geq 0$.

We now estimate the last two terms. Firstly, 

\begin{equation}
\begin{split}
&\int_{\Gamma}^{t}\lVert u_{3} \rVert_{L^{p_{1}}}^{\frac{2p_{1}}{p_{1} - 2}} X^{\frac{p_{1} - 3}{p_{1} - 2}}Z^{\frac{1}{p_{1} - 2}} d\tau\\
\lesssim& \left( \int_{\Gamma}^{t} \lVert u_{3} \rVert_{L^{p_{1}}}^{\frac{2p_{1}}{p_{1} - 3}}X d\tau \right)^{\frac{p_{1} - 3}{p_{1} - 2}}\left( \int_{\Gamma}^{t} Z d\tau \right)^{\frac{1}{p_{1} -2}}\\
\lesssim& \sup_{\tau \in [\Gamma, t]}X^{\frac{3p_{1} - 10}{4(p_{1} - 2)}}(\tau) \left(\int_{\Gamma}^{t} \lVert u_{3} \rVert_{L^{p_{1}}}^{\frac{8p_{1}}{3p_{1} - 10}} d\tau \right)^{\frac{3p_{1} - 10}{4(p_{1} - 2)}} \left(\int_{\Gamma}^{t} X d\tau \right)^{\frac{1}{4}} \left( \int_{\Gamma}^{t} Z d\tau \right)^{\frac{1}{p_{1} -2}}\\
\lesssim& \sup_{\tau \in [\Gamma, t]} X^{\frac{3}{4}}(\tau) + \left(\int_{\Gamma}^{t} Z d\tau \right)^{\frac{3}{4}}
\end{split}
\end{equation}
by H$\ddot{o}$lder's inequalities, (11), (25) and Young's inequalities. 

Next, for any $\epsilon > 0$ we estimate 

\begin{equation}
\begin{split}
& \left(\int_{\Gamma}^{t} \lVert \nabla_{h} b\rVert_{L^{2}}^{\frac{4}{3}} \lVert \nabla b\rVert_{L^{2}}^{\frac{2}{3}} \lVert u_{3} \rVert_{L^{\frac{6p_{0}}{6-p_{0}}}}^{2} d\tau   \right)^{\frac{p_{0}}{p_{0} - 2}} \int_{\Gamma}^{t}X^{\frac{p_{0} - 3}{p_{0} - 2}}Z^{\frac{1}{p_{0} - 2}} d\tau\\
\lesssim& \left(\int_{\Gamma}^{t} \lVert \nabla_{h} b\rVert_{L^{2}}^{\frac{4}{3}} \lVert \nabla b\rVert_{L^{2}}^{\frac{2}{3}} \lVert u_{3} \rVert_{L^{\frac{6p_{0}}{6-p_{0}}}}^{2} d\tau   \right)^{\frac{p_{0}}{p_{0} - 2}}  \left(\int_{\Gamma}^{t} Z d\tau \right)^{\frac{1}{p_{0} -2}}\\
\lesssim& \sup_{\tau \in [\Gamma, t]} \lVert \nabla_{h} b\rVert_{L^{2}}^{\frac{2}{3}(\frac{3p_{0} - 10}{p_{0} - 2})}\left(\int_{\Gamma}^{t}  \lVert \nabla b\rVert_{L^{2}}^{\frac{2}{3}(\frac{10}{p_{0}})} \lVert u_{3} \rVert_{L^{\frac{6p_{0}}{6-p_{0}}}}^{2} d\tau   \right)^{\frac{p_{0}}{p_{0} - 2}}  \left(\int_{\Gamma}^{t} Z d\tau \right)^{\frac{1}{p_{0} -2}}\\
\lesssim& \sup_{\tau \in [\Gamma, t]} \lVert \nabla_{h} b\rVert_{L^{2}}^{\frac{2}{3}(\frac{3p_{0} - 10}{p_{0} - 2})} \left(\int_{\Gamma}^{t} Z d\tau \right)^{\frac{1}{p_{0} -2}}\\
\leq& \epsilon \sup_{\tau \in [\Gamma, t]} \lVert \nabla_{h} b\rVert_{L^{2}}^{2} + c \left(\int_{\Gamma}^{t}Z d\tau \right)^{\frac{3}{4}}
\end{split}
\end{equation}
by H$\ddot{o}$lder's inequality, (11), (25) and Young's inequality. 

Applying (33), (34) in (32), after absorbing for $\epsilon > 0$ sufficiently small we obtain 

\begin{equation}
\sup_{\tau \in [\Gamma, t]}W (\tau) + \int_{\Gamma}^{t}Y d\tau 
\lesssim \sup_{\tau \in [\Gamma, t]} X^{\frac{3}{4}} + \left(\int_{\Gamma}^{t} Z d\tau \right)^{\frac{3}{4}}. 
\end{equation}

We are now ready to complete the $H^{1}$-bound estimate. We take $L^{2}$-inner products on (1a)-(1b) with $(-\Delta u, -\Delta b)$ to estimate 

\begin{equation*}
\begin{split}
\frac{1}{2} \partial_{t} X + Z 
\lesssim& \int (\lvert \nabla_{h}  u\rvert + \lvert \nabla_{h} b\rvert) (\lvert \nabla u\rvert^{2} + \lvert \nabla b\rvert^{2})\\
\lesssim& (\lVert \nabla_{h} u\rVert_{L^{2}} + \lVert \nabla_{h} b\rVert_{L^{2}})(\lVert \nabla u\rVert_{L^{4}}^{2} + \lVert \nabla b\rVert_{L^{4}}^{2})\\
\lesssim& W^{\frac{1}{2}} ( \lVert \nabla u\rVert_{L^{2}}^{\frac{1}{2}} \lVert \nabla u\rVert_{L^{6}}^{\frac{3}{2}} + \lVert \nabla b\rVert_{L^{2}}^{\frac{1}{2}} \lVert \nabla b\rVert_{L^{6}}^{\frac{3}{2}}) 
\lesssim W^{\frac{1}{2}} X^{\frac{1}{4}} Y^{\frac{1}{2}} Z^{\frac{1}{4}}
\end{split}
\end{equation*}
by Lemma 2.2, H$\ddot{o}$lder's and interpolation inequalities and (15). Integrating in time over $[\Gamma, \tau], \tau \in [\Gamma, t]$ and taking $\sup_{\tau \in [\Gamma, t]}$ we obtain  

\begin{equation*}
\begin{split}
&\frac{1}{2} \sup_{\tau \in [\Gamma, t]} X(\tau) + \int_{\Gamma}^{t} Z d\tau\\
\lesssim& \frac{1}{2} X(\Gamma) + \sup_{\tau \in [\Gamma, t]} W^{\frac{1}{2}}(\tau) \left(\int_{\Gamma}^{t} X d\tau \right)^{\frac{1}{4}} \left(\int_{\Gamma}^{t} Y d\tau \right)^{\frac{1}{2}} \left(\int_{\Gamma}^{t} Z d\tau\right)^{\frac{1}{4}}\\
\leq& c + c\tilde{\epsilon}^{\frac{1}{4}} (\sup_{\tau \in [\Gamma, t]} X^{\frac{3}{4}} + \left(\int_{\Gamma}^{t} Z d\tau \right)^{\frac{3}{4}} ) \left(\int_{\Gamma}^{t} Z d\tau \right)^{\frac{1}{4}}\\
\leq& c + c\tilde{\epsilon}^{\frac{1}{4}} \left( \sup_{\tau \in [\Gamma, t]} X(\tau) + \int_{\Gamma}^{t}Z d\tau \right)
\end{split}
\end{equation*} 
by H$\ddot{o}$lder's inequalities, (25), Young's inequalities and (35). Hence, for $\tilde{\epsilon}$ sufficiently small, after absorbing, we obtain 

\begin{equation*}
\sup_{\tau \in [\Gamma, t]} X(t) + \int_{\Gamma}^{t} Z d\tau \lesssim 1. 
\end{equation*}

This completes the proof of Theorem 1.1 in case $p_{1} \in (\frac{15}{2}, \infty), p_{2} \in (\frac{3}{2}, \infty)$. We now make a remark on the case $p_{2} = \frac{3}{2}$. We see that the estimates in (26), (27) both go through identically while in (28) we can estimate   

\begin{equation*}
\begin{split}
VII \lesssim& \lVert \nabla_{h} b\rVert_{L^{6}}^{2} \lVert \nabla_{h}^{2} \Delta_{h}^{-1} \omega_{3} \rVert_{L^{\frac{3}{2}}} + \lVert \nabla_{h}^{2} \Delta_{h}^{-1} u_{3} \rVert_{L^{p_{1}}} \lVert \nabla\nabla_{h} b\rVert_{L^{2}} \lVert \nabla_{h} b\rVert_{L^{\frac{2p_{1}}{p_{1} - 2}}}\\
\leq& \epsilon \lVert \nabla\nabla_{h} b\rVert_{L^{2}}^{2}+ c \lVert u_{3} \rVert_{L^{p_{1}}} \lVert \nabla\nabla_{h} b\rVert_{L^{2}} \lVert \nabla_{h} b\rVert_{L^{\frac{2p_{1}}{p_{1} - 2}}}
\end{split}
\end{equation*}
due to the homogeneous Sobolev embedding of $\dot{H}^{1}(\mathbb{R}^{3}) \hookrightarrow L^{6}(\mathbb{R}^{3})$, (29) and the smallness hypothesis. Thus, identical estimates in (30) on the second term leads to  

\begin{equation*}
\frac{1}{2} \partial_{t} W + Y 
\leq 3 \epsilon Y + c\left( 
\lVert u_{3} \rVert_{L^{p_{1}}}^{\frac{2p_{1}}{p_{1} - 2}} X^{\frac{p_{1} - 3}{p_{1} - 2}}Z^{\frac{1}{p_{1} - 2}}  + \lVert b_{3} \rVert_{L^{p_{0}}}^{\frac{2p_{0}}{p_{0} - 2}} X^{\frac{p_{0} - 3}{p_{0} - 2}} Z^{\frac{1}{p_{0} - 2}}\right).
\end{equation*}
After absorbing and integrating in time $[\Gamma, \tau]$, we obtain (31). The rest of the proof is identical.

\end{document}